\documentclass[11pt]{amsart}
\usepackage{latexsym,amsmath,amsfonts,amssymb, cite}

\title[Quadratic Weyl group multiple Dirichlet series]{Weyl group
multiple Dirichlet series constructed from quadratic characters}

\author{Gautam Chinta}
\author{Paul E. Gunnells}

\address{Department of Mathematics,
The City College of CUNY,
New York, NY 10031, USA
}

\address{Mathematisches Institut,
        Georg-August-Universit\"at,
        Bunsenstr. 3--5,
        D--37073 G\"ottingen,
        Germany}
\email{chinta@uni-math.gwdg.de}

\address{Department of Mathematics and Statistics, University of
Massachusetts, Amherst, MA 01003, USA}
\email{gunnells@math.umass.edu}

\thanks{Both authors are partially supported by the NSF.  The first
named author is also supported by a Humboldt Research Fellowship.}

\date{27 April 2006}

\DeclareMathOperator{\GL}{GL}
\DeclareMathOperator{\ord}{ord}

\DeclareMathOperator{\Supp}{Supp}
\DeclareMathOperator{\adj}{adj.}

\newcommand{\WMD}{\text{WMD}}

\newenvironment{proofof}[1]{\vspace*{.1in}\noindent{\em
Proof{#1}. \/}}{\qed\vspace{3ex}}
\begin{document}

\begin{abstract}  We construct multiple Dirichlet series in several
  complex variables whose coefficients involve quadratic residue
  symbols.  The series are shown to have an analytic continuation
  and satisfy a certain group of functional equations.  These are the
  first examples of an infinite collection of
  unstable {\it Weyl group multiple Dirichlet series} in greater
than two variables having the properties predicted in \cite{wmd1}.
\end{abstract}

\maketitle

\newcommand{\Q}{\ensuremath{{\mathbb Q}}}
\newcommand{\R}{\ensuremath{{\mathbb R}}}
\newcommand{\A}{\ensuremath{{\mathbb A}}}
\newcommand{\C}{\ensuremath{{\mathbb C}}}
\newcommand{\Z}{\ensuremath{{\mathbb Z}}}
\newcommand{\x}{\ensuremath{{\bf { x}}}}
\newcommand{\s}{\ensuremath{{\bf { s}}}}
\newcommand{\bbf}{\ensuremath{{\bf { f}}}}
\newcommand{\E}{{\mathcal E}}
\newcommand{\II}{\ensuremath{{\bf { I}}}}
\newcommand{\f}{\ensuremath{{\mathfrak { f}}}}
\newcommand{\fS}{\ensuremath{{\mathfrak { S}}}}

\newcommand{\im}{\mbox{Im}}
\newcommand{\sgn}{\mbox{sgn}}
\newcommand{\res}{\mathop{\mbox{Res}}}
\newcommand{\re}{\mbox{Re}}
\newcommand{\nin}{\noindent}
\newcommand{\mb}{\medbreak}
\newcommand{\half}{{\smfrac{1}{2}}}
\newcommand{\G}{{\gamma}}
\newcommand{\B}{{\beta}}
\newcommand{\I}{{\mathcal I}}
\newcommand{\J}{{\mathcal J}}
\newcommand{\pr}{{\prime}}
\newcommand{\norm}[1]{\left|\left|#1\right|\right|}
\newcommand{\abs}[1]{\left|#1\right|}
\newcommand{\dlim}[1]{\displaystyle {\lim_{#1}}}
\newcommand{\dsup}[1]{\displaystyle {\sup_{#1}}}
\newcommand{\dinf}[1]{\displaystyle {\inf_{#1}}}
\newcommand{\dmax}[1]{\displaystyle {\max_{#1}}}
\newcommand{\dmin}[1]{\displaystyle {\min_{#1}}}
\newcommand{\dprod}[2]{\displaystyle {\prod_{#1}^{#2}}}
\newcommand{\dcoprod}[2]{\displaystyle {\coprod_{#1}^{#2}}}
\newcommand{\dsum}[2]{\displaystyle {\sum_{#1}^{#2}\:}}
\newcommand{\prs}[2]{\left(\frac{#1}{#2}\right)}
\newcommand{\smfrac}[2]{{\textstyle \frac{#1}{#2}}} 
\newcommand{\nn}{\nonumber}
\newcommand{\no}{{\mathbb N}}
\newcommand{\cO}{{\mathcal O}}

\newtheorem{theorem}{Theorem}
\newtheorem{lemma}[theorem]{Lemma}
\newtheorem{prp}[theorem]{Proposition}
\newtheorem{corollary}[theorem]{Corollary}
\newtheorem{conjecture}[theorem]{Conjecture}
\newtheorem{fact}[theorem]{Fact}
\newtheorem{claim}[theorem]{Claim}
\theoremstyle{definition}
\newtheorem{example}[theorem]{Example}
\newtheorem{remark}[theorem]{Remark}

\newtheorem{defin}[theorem]{Definition}

\numberwithin{theorem}{section}
\numberwithin{equation}{section}

\section{Introduction}
\label{sec:intro}

Let $\Phi$ be an irreducible root system of rank $r$ with Weyl group
$W$, and let $K$ be a global field containing the $n^{th}$ roots of
unity.  In \cite {wmd1} is described a heuristic method to
associate to this data a multiple Dirichlet series $Z$ in $r$ complex
variables with coefficients given by $n^{th}$ order Gauss sums.
Moreover, $Z$ is expected to have an analytic continuation to $\C^r$
and to satisfy a group of functional equations isomorphic to $W$.
These {\it Weyl group multiple Dirichlet series} generalize several
constructions of multiple Dirichlet series that have previously
appeared in the literature.  We present some examples and outline
their connections with analytic number theory and automorphic forms at
the end of this introduction.  The paper \cite{wmd1} suggests a method
for the unified treatment of all of these examples.

Brubaker, Bump, and Friedberg \cite{wmd2} have given a precise
definition of $Z$ in the \emph{stable case}; by definition, this means
$n$ is sufficiently large for a fixed $\Phi$.  In \cite{wmd2} the
authors show that for such $n$, the Weyl group multiple Dirichlet
series admit meromorphic continuation and have the expected group of
functional equations.  They also prove that the coefficients of the
stable series satisfy a certain \emph{twisted multiplicativity}
(cf. \eqref{eq:tm}) that reduces their computation to the case of the
$P$-power coefficients, for $P$ a prime in the ring of integers of
the field $K$.  After multiplying by appropriate normalizing zeta
factors, the authors show that the $P$-parts form a Dirichlet
polynomial whose non-zero coefficients are naturally parametrized by
the elements of the Weyl group $W$.

In the unstable range, when $n$ is small relative to $\Phi$, one still
expects to be able to construct multiple Dirichlet series with the
same properties as in the stable range.  However, simple examples show
that the coefficients of the $P$-parts are no longer parametrized
solely by the elements of $W$.  We expect that terms
corresponding to the elements of $W$ will be present, but will form
only part of the $P$-part polynomial.  Some progress---including a
beautiful conjectural description of the $P$-parts for type $A$
unstable series via Gelfand--Tsetlin patterns---is given in Brubaker,
Bump, Friedberg, Hoffstein \cite{wmd3}.  There one can also find a
proof that the multiple Dirichlet series associated to $\Phi = A_2$
are the Whittaker coefficients of Eisenstein series on the metaplectic
cover of $\GL_3$.  We refer the reader to \cite{wmd3} for further
discussion of the connection between multiple Dirichlet series and
Whittaker coefficients of metaplectic Eisenstein series.

In this paper, we focus on the case $n=2$ and $\Phi$ simply-laced of
rank $r$. This assumption on $\Phi$ is no loss of generality, since
when $n=2 $ the non-simply-laced cases described in \cite{wmd1} can be
obtained by setting variables to be equal in a series associated to a
simply-laced root system.  These series are unstable for $r\geq 3$,
and the results of \cite{wmd1, wmd2} cannot be applied to derive the
desired properties of $Z$.  For such series the quadratic Gauss sums
are essentially quadratic residue symbols, and the associated Weyl
group multiple Dirichlet series can be expressed as sums of quadratic
Dirichlet $L$-functions.  Our main result is that these quadratic Weyl
group multiple Dirichlet series have a meromorphic continuation to
$\C^r$ and satisfy a group of functional equations isomorphic to the
Weyl group $W$. We refer to Section \ref{sec:mds} for the exact
definition of $Z$ and to Theorems \ref{thm:feOfZ} and \ref{thm:main}
for a precise statement of these results.

For an example, let $\Phi=A_r$ and $K=\Q$.  Then the series $Z$ has
the form
\[
\sum \frac{a(m_1,m_2,\ldots,m_r)}{m_1^{s_1}m_2^{s_2}\cdots
m_r^{s_r}},
\]
where the sum is over all positive integers $m_i$.  If $m_1m_2\cdots
m_r$ is odd and squarefree, we have
$$a(m_1,m_2,\ldots, m_r)=\prs{m_1}{m_2}\prs{m_2}{m_3}\cdots
\prs{m_{r-1}}{m_r}.$$ The coefficients satisfy the following twisted
multiplicativity property:
\begin{multline}\label{eq:tm}
a(m_1m_1^\prime,\ldots, m_rm_r^\prime)=\\
a(m_1,\ldots, m_r)a(m_1^\prime,\ldots, m_r^\prime)\prod_{j=1}^{r-1}
\prs{m_j}{m_{j+1}^\prime}\prs{m_j^\prime}{m_{j+1}},
\end{multline}
when $(m_1\cdots m_r, m_1^\prime\cdots m_r^\prime)=1$. The identity
\eqref{eq:tm} reduces the description of the coefficients to that of
the coefficients
\begin{equation}\label{eq:ppart}
a(p^{k_1}, \ldots, p^{k_r}),
\end{equation}
where $p$ ranges over all primes.  For a fixed prime $p$, the
coefficients \eqref{eq:ppart} can be organized into a generating
function
\begin{equation}\label{eq:genfun}
\sum_{k_1,\ldots,k_r\geq 0} \frac {a(p^{k_1},\ldots,
  p^{k_r})}{p^{k_1s_1}\ldots p^{k_rs_r}}.
\end{equation}
One of our main results is an explicit description of this generating
function, given in Theorem \ref{thm:invFunc}; it turns out that
\eqref{eq:genfun} is a \emph{rational} function of the $p^{-s_i}$ that
is itself invariant under a certain Weyl group action.


We conclude this introduction with a few examples of quadratic
multiple Dirichlet series that have previously appeared in the
literature.  We refer the reader to the survey papers \cite{cfh2,
wmd1} for a more comprehensive list of examples.

The first example (in more
than one variable) was found by Siegel \cite{si}:
\begin{equation}\label{eqn:siegelZ}
Z(s,w)=
\mathop{\sum\sum}_{\substack{d,m\geq 1\\ d, m\ \text{odd}}}
\frac{\chi_d(\hat
  m)}{m^sd^w}a(d,m),
\end{equation}
where $\hat m$ denotes the part of $m$ relatively prime to the
squarefree part of $d$ and $\chi_d$ is the quadratic character
associated to the extension $\Q(\sqrt{d})$ of $\Q$.    The
multiplicative factor $a(d,m)$ is defined by
$$a(d,m)=\prod_{\substack{\text{$p$ prime}\\p^k||d,\ p^l||m}}
a(p^k, p^l),$$
and
\begin{equation}\label{eqn:siegel}
a(p^k,p^l)=\left\{ \begin{array}{ll}
\min(p^{k/2}, p^{l/2}) & \mbox{if $\min(k, l)$ is even,} \\
0 & \mbox{otherwise.}
\end{array}\right.
\end{equation}
Siegel obtained this series as the Mellin transform of a half-integral
weight Eisenstein series for the congruence subgroup $\Gamma_0(4)$.
(Actually, Siegel's series is a linear combination of series of this
form.)  As Siegel noted, this integral representation implies two
functional equations for $Z(s,w)$, one coming from the functional
equation of the Eisenstein series, and one coming from the Mellin
transform, via the automorphicity of the Eisenstein series.  These
functional equations take the form
\[
Z(s,w)\mapsto Z(1-s, w+s-1/2) \quad \text{and}\quad Z(s,w)\mapsto
Z(s, 3/2-s-w).
\]
These two functional equations commute with one
another, and thus generate a group isomorphic to the Klein four group.

In fact, it turns out that Siegel's series satisfies a group of
\emph{twelve} functional equations.  In our notation, Siegel's series
is the quadratic series associated to the root system of type $A_2$.
This means that \eqref{eqn:siegelZ} actually possesses a group of
functional equations $G$ isomorphic to the direct product of the Weyl
group of type $A_{2}$ together with order 2 group of symmetries of the
Dynkin diagram of $A_2.$ These extra functional equations---which are
not at all obvious from Siegel's presentation of his series---were
first noted in unpublished work of Bump and Hoffstein, who recognized
this Mellin transform of the metaplectic $GL_2$ Eisenstein series as
the Fourier--Whittaker coefficient of a minimal parabolic metaplectic
Eisenstein series on the double cover of $GL_3$.\footnote{A more
general connection between double Dirichlet series and Whittaker
coefficients of a metaplectic $GL_3$ Eisenstein series is proven in
\cite{wmd3}.}  The full group of functional equations, as well as the
meromorphic continuation of $Z(s,w)$, was worked out in detail by
Fisher--Friedberg \cite{ff}, using methods totally separate from the
work of Bump--Hoffstein.  For an application of $Z(s,w)$ to a mean
value result for sums of quadratic Dirichlet $L$-functions, see
Goldfeld--Hoffstein \cite{gh} as well as \cite{ff}.

For a rank $3$ example, take the Rankin--Selberg convolution
of two half-integral weight Eisenstein series for $\Gamma_0(4)$.  This
yields the quadratic $A_3$ series, which has the form
\begin{equation}\label{eqn:a3Z}
\mathop{\sum\sum\sum}_{\substack{d,n_1,n_2>0\\
d, n, n_{2}\ \text{odd}}} \frac{\chi_d(\hat
n_1)\chi_d(\hat n_2)}{n_1^{s_1}n_2^{s_2}d^w}a(n_1,n_2,d).
\end{equation}
Here $a(n_1, n_2, d)$ is a multiplicative weighting factor first
explicitly written down by Fisher and Friedberg \cite{ff2} (see also
our Example \ref{ex:f}).  It is expected that the $A_3$ series is a
Whittaker coefficient of a minimal parabolic Eisenstein series on the
double cover of $GL_4$.  Applications of the $A_3$ series include mean
value results for sums of squares of quadratic Dirichlet
$L$-functions.

More examples of higher rank have also appeared in the literature and
have been applied to analytic number theory.  The quadratic $D_4$
series was treated by Diaconu, Goldfeld and Hoffstein \cite{dgh}, who
used it to prove mean value results for sums of cubes of quadratic
Dirichlet $L$-functions.  This was first proved by Soundararajan
\cite{so} by other methods.  The results of \cite{dgh} and \cite{so}
are stated over \Q, but the methods of multiple Dirichlet series work
over any global field.  One of us (GC) recently used the quadratic
$A_5$ series to establish a mean value result for central values of
zeta functions of biquadratic number fields \cite{chbq}.  The results
of this paper simultaneously unify and generalize all of these earlier
constructions.

Finally, we remark that one may also construct double Dirichlet series
roughly of the form
\begin{equation}\label{eqn:zswf}
\mathop{\sum\sum}_{\substack{d,n\geq 0\\
d, n\ \text{odd}}} \frac{\chi_d(\hat
n)}{n^sd^w}b_g(n),
\end{equation}
where the $b_g(n)$ are Fourier coefficients of a Hecke cuspform $g$ on
$GL_2$ or $GL_3$. These have been studied in the papers \cite{ghp,
bfh90} (for $g$ on $GL_2$), and \cite{bfh04,cd} (for $g$ on $GL_3$).
Though we do not directly address such series in this paper, our
methods may easily be adapted to establish the analytic continuation
and functional equations of \eqref{eqn:zswf}.

We briefly indicate how to define \eqref{eqn:zswf} precisely when $g$
is a $GL_3$ form.  This is the heart of the problem, since once the
series has been correctly defined, it is easy to mimic the procedure
of Section 5 to establish the functional equations and analytic
continuation.  To precisely define \eqref{eqn:zswf}, it once again
suffices to specify its $p$-part.  Let $\alpha_1,\alpha_2,\alpha_3$ be
the Satake parameters of $g$ at an unramified prime $p$. Let
$f(x_1,x_2,x_3,x_4)$ be the rational function associated to the root
system of type $D_4$ given by Theorem \ref{thm:invFunc}.  (We take
$x_4$ to be the variable corresponding to the central node.)  Then the
generating series giving the precise form of the $p$-part of the
series \eqref{eqn:zswf} is $f(\alpha_1 x,\alpha_2 x,\alpha_3 x, y)$,
where $x=p^{-s}, y=p^{-w}$.

\medskip \noindent \textbf{Acknowledgments}.  We thank Jim Humphreys
for helpful conversations.  We are grateful to Ben Brubaker, Dan Bump,
Sol Friedberg, and Jeff Hoffstein for making available to us their
preprint \cite{wmd3} and for enlightening correspondence.  We also
thank the organizers of the Bretton Woods Workshop on Multiple
Dirichlet Series (July 2005), where some of this work was carried out.
Finally we thank an anonymous referee for many helpful clarifications.

\section{Preliminaries} \label{sec:prelims}
Let $K$ be a number field with ring of integers $\cO$. Let $S_f$ be a
finite set of non-archimedean places such that $S_f$ contains all
places dividing $2$ and the ring of $S_f$-integers $\cO_{S_f}$ has
class number $1$. Let $S_\infty$ be the set of archimedean places, and
let $S=S_f\cup S_\infty$.

Let $\prs{a}{*}$ be the quadratic residue symbol attached to the
extension $K(\sqrt{a})$ of $K$, extended as in
\cite{ff}; we review the definition below.  A slightly
different but essentially equivalent formalism appears in the
papers \cite{wmd1, wmd2, wmd3}.  We find the setup of
\cite{ff} simpler in the quadratic case.

For each place $v$, let $K_v$ denote the completion of $K$ at $v$.
For $v$ nonarchimedean, let $P_v$ be the corresponding ideal of $\cO$,
and let $q_v = |P_v|$ be its norm.  Let $C$ be the product
$\prod_{v\in S_f}P_v^{n_v}$ where $n_v$ is defined to be
$\max\{\mbox{ord}_v(4),1\}$. Let $H_C$ be the narrow ray class group
modulo $C$, and let $R_C=H_C\otimes \Z/2\Z$.  Write the finite group
$R_C$ as a direct product of cyclic groups, choose a generator for
each, and let $\E_0$ be a set of ideals of $\cO$ prime to $S$ that
represent these generators.  For each $E_0\in \E_0$, choose
$m_{E_0}\in K^\times$ such that $E_0\cO_{S_f}=m_{E_0}\cO_{S_f}$.  Let
$\E$ be a full set of representatives for $R_C$ of the form
$\prod_{E_0\in\E_0}E_0^{n_{E_0}}$, with $n_{E_0}\in\Z$.  If
$E=\prod_{E_0\in\E_0}E_0^{n_{E_0}}$ is such a representative, then let
$m_E=\prod_{E_0\in\E_0}m_{E_0}^{n_{E_0}}$.  Note that
$E\cO_{S_f}=m_E\cO_{S_f}$ for all $E\in\E$.  For convenience we
assume that $\cO\in \E$ and $m_\cO=1$.

Let $\J(S)$ be the group of fractional ideals of $\cO$ coprime to
$S_f$.  Let $I,J\in \J(S)$ be coprime.  Write $I=(m)EG^2$ with
$E\in\E$, $m\in K^\times$, $m\equiv 1\bmod C$, and $G\in \J(S)$ such
that $(G,J)=1$.  Then, following \cite{ff}, the quadratic residue
symbol $\prs{m m_E}{ J}$ is defined, and if $ I = (m')E'{G'}^2 $ is
another such decomposition, then $E'=E$ and $ \prs{m' m_E}{ J} =
\prs{m m_E}{ J}$.  In view of this define the quadratic residue symbol
$\prs{I}{J}$ to be $\prs{m m_E}{J}$. For $I=I_0I_1^2$ with $I_0$
squarefree, we denote by $\chi_I$ the character
$\chi_I(J)=\chi_{I_0}(J)=\prs{I_0}{J}$.  Further, in the expression
$\chi_I(\hat J)$, we let $\hat J$ represent the part of $J$ coprime to
$I_0$. This character $\chi_I$ depends on the choices above, but we
suppress this from the notation.

%
\begin{prp}[Reciprocity]\cite{ff} Let $I$,~$J \in \J(S) $ be coprime, and $
\alpha(I,J) = \chi_I(J)\chi_{J}(I)^{-1}$.  Then $\alpha(I,J)$ depends
only on the images of $I$~and~$J$ in $ R_C $.
\end{prp}

\begin{proof} See Neukirch \cite{ne}, Theorem 8.3 of Chapter 6.
\end{proof}

Let $\I(S)$ be the set of integral ideals prime to
$S_f$. Let $L^S(s,\chi_J)$ be the $L$-function of the character
$\chi_J$, with the places in $S$ removed.  We let $L_S(s,\chi_J)$ be
the product over the places in $S$.  Thus
$$L(s,\chi_J)=L^S(s,\chi_J)L_S(s,\chi_J).$$

If $\xi$ is any id\`ele class character then the
completed $L$-function $L(s,\xi)$ satisfies a functional
equation
%
%
\begin{equation}
\label{eq:LFE}
L(s,\xi)=\epsilon(s,\xi)L(1-s,\xi^{-1}),
\end{equation}
where $\epsilon(s,\xi)$ is the epsilon factor of $\xi$.
%
%
\begin{prp}
\label{prop:epsilonFE}
Let $E,J\in \cO(S)$ be squarefree with
associated characters $\chi_E,\chi_J$ of conductors $\mathfrak
f_E,\mathfrak f_J$ respectively. Suppose that $\chi_J=\chi_E\chi_I$
with $I\in K^\times$, $I\equiv 1 \bmod C$. Let $\psi$ be another
character unramified outside $S$. Then
\begin{equation}
\epsilon(s,\chi_J\psi)=\epsilon(1/2,\chi_I)\psi(\mathfrak f_J/
  \mathfrak f_E)\left(|\mathfrak f_J/
  \mathfrak f_E| \right) ^{1/2-s} \epsilon(s,\chi_E\psi).
\end{equation}
\end{prp}

Here $\epsilon(1/2,\chi_I)$ is given by a (normalized) Gauss sum, as
in Tate's thesis.  When $\chi_I$ is a quadratic character, we have
$\epsilon(1/2, \chi_I)=1$.

We remark that the $\Gamma$-factors of the $L$-function appear in
the contribution of the archimedean places $L(s,\chi_J).$  When
the base field $K$ is totally real, these $\Gamma$-factors will
depend on $\chi_J,$ but only on the narrow ray class of $J.$  For
example, when $K=\Q$ and $d$ a fundamental discriminant, the
$\Gamma$-factor of $L(s,\chi_d)$ is $\Gamma(\frac s2)$ if $d>0$
and $\Gamma(\frac{s+1}{2})$ if $d<0.$

\section{A Weyl group action on rational
  functions}\label{sec:rationalFunctions}
Let $\Phi$ be an irreducible simply laced root system of rank $r$ with
Weyl group $W$. Choose an ordering of the roots and let $\Phi = \Phi^+
\cup \Phi^{-}$ be the decomposition into positive and negative roots.
 Let
$$\Sigma=\{\alpha_1, \alpha_2, \ldots, \alpha_r\}$$ be the set of
simple roots and let $\sigma_i$ be the Weyl group element
corresponding to the reflection through the hyperplane perpendicular
to $\alpha_i$.  We say that $i$ and $j$ are \emph{adjacent} if $i\neq
j$ and $(\sigma_i \sigma_j)^3=1$.  The Weyl group $W$ is generated by
the simple reflections $\sigma_1,\sigma_2,\ldots, \sigma_r$, which
satisfy the relations
\begin{equation}\label{eqn:WRelations}
(\sigma_i\sigma_j)^{r(i,j)} = 1
\mbox {\ with\ } r(i,j)=\left\{\begin{array}{ll}
3 & \mbox{\ if $i$ and $j$ are adjacent,}\\
1 & \mbox{\ if $i=j$, and}\\
2 & \mbox{\ otherwise,}
\end{array}\right.
\end{equation}
for $1\leq i, j\leq r$.
The action of the generators $\sigma_i$ on the roots is
\begin{equation}\label{eqn:WactionRoots}
\sigma_i\alpha_j=\left\{\begin{array}{ll}
\alpha_i+\alpha_j & \mbox{\ if $i$ and $j$ are adjacent,}\\
-\alpha_j & \mbox{\ if $i=j$, and  }\\
\alpha_j & \mbox{\ otherwise.}
\end{array}\right.
\end{equation}
Though it will play no role in this section, we will assume that the
indices are ordered so that for each $j$, the $i$ adjacent to $j$ are
either all less than $j$ or all greater than $j$.

Let $l$ denote the length function on $W$ with respect to the
generators $\sigma_1, \sigma_2,\ldots, \sigma_r$, and define
$$\sgn(w)=(-1)^{l(w)}.$$
Let $\Lambda_{\Phi}$ be the lattice generated by the roots.  Any
$\alpha \in \Lambda_{\Phi}$ has a unique representation as an
integral linear combination of the simple roots:
\begin{equation}\label{eqn:repn}
\alpha=k_1\alpha_1+k_2\alpha_2+\cdots+k_r\alpha_r.
\end{equation}
We call the set $\Supp (\alpha)$ of $j$ such that $k_{j} \not = 0$ in
\eqref{eqn:repn} the \emph{support} of $\alpha$.
We put
$$d(\alpha)=k_1+k_2+\cdots+k_r.$$
Introduce a partial ordering on $\Lambda_{\Phi}$ by defining
$\alpha\succeq 0$ if each $k_i\geq 0$ in \eqref{eqn:repn}.  Given
$\alpha, \beta\in \Lambda_{\Phi}$, define
$\alpha\succeq\beta$ if $\alpha-\beta\succeq 0$.

Let
$$\rho=\frac 12\sum_{\alpha\in \Phi^+}\alpha$$ be half the sum of the
positive roots.
For each $w$ in the Weyl group set
$$\Phi(w)=\{\alpha\in \Phi^+:w(\alpha)\in \Phi^-\}.$$
We gather some simple properties of $W$ we will need later.

\begin{lemma}\label{lemma:Phi}
Let $w\in W$.
\begin{enumerate}
\item[(a)] The cardinality of $\Phi(w)$ is the length $l(w)$ of $w$.
\item[(b)] We have
\begin{equation}
\rho-w\rho=\sum_{\alpha\in \Phi(w^{-1})} \alpha.
\end{equation}
\item[(c)]  Let $\sigma_i\in W$ be a generator such that
  $l(\sigma_iw)=l(w)+1. $
  Then $$\Phi(\sigma_iw)=\Phi(w)\cup \{w^{-1}\alpha_i \}.$$
\item[(d)] Let $\sigma_i\in W$ be a generator such that
$l(w\sigma_i)=l(w)+1. $
  Then $$\Phi(w\sigma_i)=\sigma_i\left(\Phi(w)\right)\cup \{\alpha_i
  \}.$$
\item [(e)] The set of simple reflections $\sigma_{i}$ occurring in any
reduced expression for $w$ is uniquely determined by $w$.
\item[(f)] Let $J = \Supp (\rho -w\rho)$.  Then $w$ lies in the
subgroup $\langle \sigma_j \mid j\in J\rangle $.
\end{enumerate}
\end{lemma}

\begin{proof}
Statements (a)--(e) can be easily found in many standard references,
e.g. \cite{humph}.  We were unable to locate a precise reference
for (f), and so for the convenience of the reader provide a proof.  We
will prove (f) under the assumption that $W$ is a simply-laced Weyl
group.  Note that by (b) the set $\Supp (\rho -w\rho)$ makes sense for
any $w\in W$.   We proceed by induction on $l (w)$.

First assume $l (w)= 1$, so that $w=\sigma_{i}$, a simple reflection.
Then $\rho -\sigma_{i} \rho = \alpha_{i}$.  Hence the result is true
in this case.

Now assume the result is true for lengths up to $l (w)$.  By (e) it
suffices to check the truth of the statement on any reduced expression
for $w$.  Let $ \sigma_{i}u$ be a reduced expression for $w$, so that
$l (u) = l (w)-1$.  Let $J = \Supp (\rho -u\rho)$ and write $\rho
-u\rho = \sum_{j\in J} k_{j}\alpha_{j}$, where $k_{j}>0$.  By (b) and
(d), we have
\begin{equation}\label{eq:rhoexp}
\rho -w\rho = \sum_{\alpha \in \Phi (w^{-1})} \alpha = \alpha_{i} +
\sum_{\alpha \in \sigma_{i}\Phi (u^{-1})} \alpha = \alpha_{i}+
\sum_{j\in J}k_{j}\sigma_{i}(\alpha_{j}).
\end{equation}
Write the last expression as
\begin{equation}\label{eq:rhoexp2}
\alpha_{i}+\sum_{j=1,\dots ,r}k'_{j}\alpha_j=
\sum_{j=1,\dots ,r}k''_{j}\alpha_j.
\end{equation}
We claim $k''_{j} \not = 0$ if $j\in J$.  Indeed, assume $j\not =i$.
If $j$ is not adjacent to $i$, then $\sigma_{i} (\alpha_{j}) =
\alpha_{j}$.  On the other hand if $j$ is adjacent to $i$, then
$\sigma_{i} (\alpha_{j}) = \alpha_{i}+\alpha_{j}$.  Hence if $j\not
=i$ we must have $k''_j=k'_{j}\geq k_{j}$.

Now supppose $j=i$.  Then the only problem is that we might have
$k_{i}'=-1$, which would lead to an expression for $\rho -w\rho$ not
involving $\alpha_{i}$.  However, by \eqref{eq:rhoexp} and
\eqref{eq:rhoexp2} we have
\begin{equation}\label{eq:possum}
\sum_{\alpha \in \sigma_{i}\Phi (u^{-1})} \alpha =\sum_{j=1,\dots
,r}k'_{j}\alpha_j,
\end{equation}
and the left of \eqref{eq:possum} is a sum over positive roots.  Thus
$k''_i=1+k'_i> 0$.  This completes the proof.
\end{proof}

Let $F=\C({\bf x})=\C(x_1,x_2,\ldots, x_r)$ be the field of rational
functions in the variables $x_1,x_2,\ldots, x_r$. For any $\alpha \in
\Lambda_{\Phi}$, let $\x^\alpha\in F$ be the monomial
$x_1^{k_1}x_2^{k_2}\cdots x_r^{k_r}$, where the exponents $k_{i}$ are
determined as in \eqref{eqn:repn}.  Our immediate goal is to define an
action of the Weyl group $W$ on $F$.  It will turn out that to
construct a multiple Dirichlet series with group of functional
equations isomorphic to the group $W$, it suffices to construct a
rational function $f$ invariant under this $W$-action and satisfying
certain limiting conditions, see Section \ref{sec:fes} and Proposition
\ref{prp:LFunction}.

We define this $W$-action in stages.  First, for $\x=(x_1,x_2,\ldots,
x_r)$ define $\sigma_i\x=\x^\prime$, where
\begin{equation}\label{eqn:wiaction1}
x_j^\prime=\left\{\begin{array}{ll}
x_i x_j\sqrt{q} & \mbox{\ if $i$ and $j$ are adjacent,}\\
1/(q x_j) & \mbox{\ if $i=j$, and}\\
x_j & \mbox{\ otherwise.}
\end{array}\right.
\end{equation}
It is easy to see that 
\begin{equation}\label{eqn:sigmaandx}
\begin{array}{ll}
\sigma_{i}^{2}\x = \x&\quad \text{for all $i$,}\\
\sigma_{i}\sigma_j\sigma_{i}\x=\sigma_{j}\sigma_{i}\sigma_{j}\x &\quad \text{if
$i$ and $j$ are adjacent,}\\
\sigma_{i}\sigma_{j}\x =
\sigma_{j}\sigma_{i}\x&\quad \text{otherwise}.
\end{array}
\end{equation}
Next, define $\epsilon_i \x=\x^\prime$, where
\begin{equation}
x_j^\prime=\left\{\begin{array}{ll}
-x_j & \mbox{\ if $i$ and $j$ are adjacent,}\\
x_j & \mbox{\ otherwise.}
\end{array}\right.
\end{equation}
Clearly $\epsilon_i^{2}\x = \x$ and $\epsilon_i \epsilon_j\x =
\epsilon_j \epsilon_i \x $, and we have 
\begin{equation}\label{eqn:epsandsig}
\sigma_{i}\epsilon_{j}\x  = \left\{\begin{array}{ll}
\epsilon_{i}\epsilon_{j}\sigma_{i}\x & \mbox{\ if $i$ and $j$ are adjacent,}\\
\epsilon_{j}\sigma_{i}\x& \mbox{\ otherwise.}
\end{array}\right.
\end{equation}
For $f\in F$ define 
\begin{equation}
f_i^+(\x)=\frac{f(\x)+f(\epsilon_i \x)}{2}\mbox{ \ \ and\ \ }
 f_i^-(\x)=\frac{f(\x)-f(\epsilon_i \x)}{2}.
\end{equation}
Finally we can define the action of $W$ on $F$ for a generator
$\sigma_i\in W:$
\begin{equation}\label{wiaction}
(f|\sigma_i)(\x)=-\frac{1-qx_i}{qx_i(1-x_i)}f_i^+(\sigma_i\x)+
\frac{1}{x_i\sqrt{q}}f_i^-(\sigma_i\x)
\end{equation}

\begin{lemma}
The definition \eqref{wiaction} extends to give an action of $W$ on $F$.
\end{lemma}

\begin{proof}
The proof amounts to verifying that the relations
\eqref{eqn:WRelations} are respected by \eqref{wiaction}.  These are
straightforward computations that involve identities in rational
functions that are independent of $f$ and the global structure of the
root system $\Phi$.  We will show in detail that
$f|\sigma_{i}^{2} = f$, and will explain what computations are
involved in proving $f|\sigma_{i}\sigma_{j}\sigma_{i} =
f|\sigma_{j}\sigma_{i}\sigma_{j}$ when $i$ is adjacent to $j$.  The
final relation, that $f|\sigma_{i}\sigma_{j} = f|\sigma_{j}\sigma_{i}$
when $i$ and $j$ are not adjacent, is proved by the same technique and
will be left to the reader.

Define 
\begin{eqnarray*}
c_{i}(\bf x) &=& \frac{1}{2}\left( \frac{q x_i-1}{qx_i (1-x_i)
   }+\frac{1}{\sqrt{q} x_i}\right), \mbox{\ \ \ and} \\
d_{i}(\bf x) &=& \frac{1}{2}\left( \frac{q x_i-1}{qx_i (1-x_i)
   }-\frac{1}{\sqrt{q} x_i}\right)
\end{eqnarray*}
for $i=1,2,\ldots , r.$  We can rewrite \eqref{wiaction} as
\begin{equation}\label{eqn:Wactionrewritten}
(f|\sigma_i)(\x)=c_i({\bf x}) f(\sigma_{i}\x)+d_i({\bf x})f(\epsilon_i \sigma_{i}\x).
\end{equation}
It is then easy to compute 
\begin{multline*}
(f|\sigma_{i}^{2}) (\x) = (c_{i} (\x)c_{i} (\sigma_{i}\x) + d_{i}
(\x)d_{i} (\sigma_{i}\x)) f (\x) \\
+ (c_{i} (\x)d_{i} (\sigma_{i}\x) + d_{i}
(\x)c_{i} (\sigma_{i}\x)) f (\epsilon_{i}\x).
\end{multline*}
Hence for $f|\sigma_{i}^{2} = f$ we need 
\begin{subequations}
\begin{align}
c_{i} (\x)c_{i} (\sigma_{i}\x) + d_{i}
(\x)d_{i} (\sigma_{i}\x) &= 1,\label{eqn:conditions1}\\
c_{i} (\x)d_{i} (\sigma_{i}\x) + d_{i}
(\x)c_{i} (\sigma_{i}\x) &= 0.\label{eqn:conditions2}
\end{align}
\end{subequations}

This is quickly seen as follows.  Let 
\[
A = \frac{q x_i-1}{qx_i (1-x_i)}, \quad B = \frac{1}{\sqrt{q} x_i}.
\]
Then $c_{i} (\x) = (A+B)/2$ and $d_{i} (\x) = (A-B)/2$.  One can check
that $c_{i} (\sigma_{i}\x) = (A^{-1}+B^{-1})/2$ and $d_{i}
(\sigma_{i}\x) = (A^{-1}-B^{-1})/2$, so
\eqref{eqn:conditions1}--\eqref{eqn:conditions2} follow easily.

Now we suppose $i$ is adjacent to $j$, and consider $F_{1} =
f|\sigma_{i}\sigma_{j}\sigma_{i}$ and $F_{2}=
f|\sigma_{j}\sigma_{i}\sigma_{j}$.
Repeatedly applying \eqref{eqn:Wactionrewritten} and the relations
\eqref{eqn:sigmaandx} and \eqref{eqn:epsandsig}, we can write both $F_{1}$ and $F_{2}$ as linear
combinations of the four functions 
\[
\text{$f (\sigma_{i}\sigma_{j}\sigma_{i}\x )$, $f (\epsilon_{i}
\sigma_{i}\sigma_{j}\sigma_{i}\x )$, $f (\epsilon_{j}
\sigma_{i}\sigma_{j}\sigma_{i}\x )$, and $f (\epsilon_{i}
\epsilon_{j}\sigma_{i}\sigma_{j}\sigma_{i}\x )$.}
\]
Comparing coefficients of these linear combinations gives four
identities in rational functions that must be satisfied for $F_{1}$ to
equal $F_{2}$.  For instance, the identity needed for equality of the coefficients
of $f (\sigma_{i}\sigma_{j}\sigma_{i}\x)$ in $F_{1}, F_{2}$ is 
\begin{multline*}
c_{i} (\x)c_{j} (\sigma_{i}\x)c_{i} (\sigma_{j}\sigma_{i}\x) + d_{i}
(\x)d_{2} (\epsilon_{i}\sigma_{i}\x)d_{i}
(\epsilon_{i}\sigma_{j}\sigma_{i}\x) = \\
c_{j} (\x)c_{i} (\sigma_{j}\x)c_{j} (\sigma_{i}\sigma_{j}\x) + d_{j}
(\x)d_{2} (\epsilon_{j}\sigma_{j}\x)d_{j}
(\epsilon_{j}\sigma_{i}\sigma_{j}\x).
\end{multline*}
Such identities are easily verified with the aid of a computer algebra system.
This completes the proof.
\end{proof}

\begin{lemma}\label{lemma:Waction}
Let $g,h\in F$ and $w\in W$.
\begin{enumerate}
\item [(a)]  $(g+h)|w=g|w+h|w $
\item[(b)] If $g(\x)=g_\alpha(\x)=\x^\alpha$ is a monomial, then
$$g(w\x)=q^{ d(w\alpha-\alpha)/2}\x^{w\alpha}.$$
\item[(c)]
If $g$ is an even function of all the $x_j$, then
$$(gh|w)(\x)=g(w\x)\cdot (h|w)(\x).$$
\end{enumerate}
\end{lemma}

\begin{proof}
Each part of the Lemma can be proven by first establishing the result
for the generators $\sigma_{i}$, and then verifying that if the result
is true for $w_1,w_{2}\in W,$ then it is true for the product $w_1w_2$.
 Part
(a) is obvious.  For part (b), we have $g(\sigma_i
\x)= q^{d(\sigma_i\alpha-\alpha)/2}\x^{\sigma_i\alpha}$ by
\eqref{eqn:WactionRoots} and \eqref{eqn:wiaction1}.  Assume (b)
holds for $w_1, w_2\in W$. Then we have
\begin{eqnarray*}
w_1 w_2 (\x^\alpha)&=&q^{
  d(w_2\alpha-\alpha)/2}w_1(\x^{w_2\alpha})\\
&=&q^{ d(w_2\alpha-\alpha)/2}q^{d(w_1w_2\alpha-
w_2\alpha)/2}\x^{w_1w_2\alpha}\\
&=&q^{d((w_1w_2)\alpha-\alpha)/2}\x^{(w_1w_2)\alpha},
\end{eqnarray*}
as required.  For part (c), first note that
if $g$ is an even function of $x_j$ for each index $j$
adjacent to $i$, then
$$(gh|\sigma_i)(\x)=g(\sigma_i\x)\cdot(h|\sigma_i)(\x).$$ Part (b)
implies that if $g (\x)$ is even in any variable, then $g(\sigma_i\x)$
is even in the same variable.
The proof of (c) is now easily completed.
\end{proof}

We now state the main result of this section.

\begin{theorem}\label{thm:invFunc}
There exists a rational function $f \in F$ that is
  $W$-invariant under the $|$ operation induced by \eqref{wiaction}
and satisfies
\begin{enumerate}
\item for each $i=1,2,\ldots, r$, the function $f$  satisfies the following
  limiting condition: if $x_j= 0$ for every $j$ adjacent to $i$, then
\begin{equation}\label{limCond}
f(\x)(1-x_i) \mbox{\ is independent of \ } x_i.
\end{equation}
\item $f(0,0,\ldots, 0)=1.$
\end{enumerate}
\end{theorem}

\begin{remark}  We expect that the rational function satisfying the
  conditions of Theorem \ref{thm:invFunc} is unique.  We have verified the
  uniqueness by a laborious induction for the root systems $A_n\ (n\leq
5)$ and $D_4$.
\end{remark}

\begin{example}\label{ex:fa2}
For the root system $A_2$, the rational function $f$ satisfying the
conditions of Theorem \ref{thm:invFunc} is
\begin{equation}\label{eqn:a2series}
f_{A_{2}} = f_{A_{2}}(x_1,x_2)=\frac{1-x_1x_2}{(1-x_1)(1-x_2)(1-q x_1^2
x_2^2)}.
\end{equation}
The Taylor series coefficients of $f_{A_{2}}$ coincide with the $q$-part of
Siegel's series \eqref{eqn:siegelZ}.  That is, if we write
$$f(x_1, x_2)=\sum_{k,l\geq 0} a_{kl} (q) x_1^k x_2^l,$$
then
$$a_{kl} (q) =\left\{ \begin{array}{ll}
\min(q^{k/2}, q^{l/2}) & \mbox{if $\min(k, l)$ is even,} \\
0 & \mbox{otherwise.}
\end{array}\right.$$
This should be compared with \eqref{eqn:siegel}.
\end{example}

\begin{example}\label{ex:f}
For the root system $A_3$, with central node corresponding to
  $x_2$, the rational function is
$f_{A_{3}} = f_{A_{3}}(x_1,x_2, x_3)=$
$$\frac{1-x_1x_2-x_2x_3+x_1x_2x_3+
qx_1x_2^2x_3-qx_1^2x_2^2x_3-qx_1x_2^2x_3^2+qx_1^2x_2^3x_3^2}
{(1-x_1)(1-x_2)(1-x_3)(1-q x_1^2
  x_2^2)(1-qx_2^2x_3^2)(1-q^2x_1^2x_2^2x_3^2)}
.$$
This can be expressed in terms of the $A_2$ rational function
$f_{A_{2}}$ from Example \ref{ex:fa2}.  Indeed, for $|x_i|<1/q$, we have
\begin{equation}\label{eqn:conv}
f_{A_{3}}(x_1,x_2,x_3)=\frac{1}{1-qx_1x_2^2x_3}\int
f_{A_2}(x_1,t)f_{A_2}(x_2t^{-1}, x_3) \frac{dt}{t},
\end{equation}
where the integral is taken over the circle $|t|=1/q$.

The identity \eqref{eqn:conv} originates in the representation of the
$A_3$ multiple Dirichlet series \eqref{eqn:a3Z} as a Rankin-Selberg
convolution of two metaplectic Eisenstein series on the double cover
of $GL_2$ (cf. Section \ref{sec:intro}).  The factor
$(1-qx_1x_2^2x_3)^{-1}$ can be interpreted as the $q$-part of the
normalizing zeta factor arising in the convolution, cf. \cite[Section
1.1]{bu}.
\end{example}

The relation between the above examples and the results of \cite{wmd1,
  wmd2, wmd3} is discussed in Remark \ref{relation-to-wmd} at the end of
  this section.

Since $W$ is finite, it is easy to construct functions in $F$ that are
$W$-invariant by averaging over the group.  The difficulty lies
in finding the proper function to average so that the condition
\eqref{limCond} is satisfied.

To this end,
define
$$\Delta(\x)=\prod_{\alpha\in \Phi^+}(1-q^{d(\alpha)}\x^{2\alpha}),$$
and let
$$j(w,\x)=\Delta(\x)/\Delta(w\x).$$
Then $j$ satisfies the one-cocycle relation
\begin{equation}\label{eqn:cocycle}
j(ww^\prime, \x)=j(w,w^\prime\x)j(w^\prime,\x).
\end{equation}

\begin{lemma}\label{lemma:j} We have
$$j(\sigma_i, \x)=-q x_i^2$$ for each simple reflection $\sigma_i$.
Moreover, let $w\in W$ and let $\alpha = \rho -w^{-1}\rho$.  Then we have
$$j(w,\x)=\sgn(w)q^{d(\alpha)}\x^{2\alpha}.$$
\end{lemma}

\begin{proof}
The second statement follows from the first and
the cocycle relation \eqref{eqn:cocycle}.  For the first, write
$$\Delta(\x)=(1-q^{d(\alpha_i)}\x^{2\alpha_i})\prod_{\substack{\alpha\in
  \Phi^+\\ \alpha\neq\alpha_i}}(1-q^{d(\alpha)}\x^{2\alpha}).$$

Using Lemma \ref{lemma:Waction},
\begin{eqnarray*}
\Delta(\sigma_i\x)&=&
(1-q^{d(\alpha_i)}q^{d(\sigma_i\alpha_i-\alpha_i)}\x^{2\sigma_i\alpha_i})
\prod_{\substack{\alpha\in \Phi^+\\ \alpha\neq\alpha_i}}
(1-q^{d(\alpha)}q^{d(\sigma_i\alpha-\alpha)}\x^{2\sigma_i\alpha})\\
&=&(1-q^{-d(\alpha_i)}\x^{-2\alpha_i})
\prod_{\substack{\alpha\in \Phi^+\\ \alpha\neq\alpha_i}}
(1-q^{d(\sigma_i\alpha)}\x^{2\sigma_i\alpha})
\end{eqnarray*}
since $\sigma_i\alpha_i=-\alpha_i$.  But by Lemma \ref{lemma:Phi}  the
positive roots in $\Phi^+\backslash \{\alpha_i\}$ are permuted by
$\sigma_i$. Therefore
$$\Delta(\sigma_i\x)=-\frac{1}{qx_i^2} \Delta(\x),$$
as claimed.
\end{proof}

We are now ready to construct the function whose existence is claimed
in Theorem \ref{thm:invFunc}.  Define
\begin{equation}\label{defoff0}
f_0(\x)=\sum_{w\in W} j(w,\x) (1|w)(\x),
\end{equation}
and put
\begin{equation}\label{eqn:defoff}
f(\x)=f_0(\x)\Delta(\x)^{-1}.
\end{equation}
We claim $f (\x)$ satisfies the conditions of Theorem
\ref{thm:invFunc}.

The invariance of $f$ is clear.  To verify the limiting condition
\eqref{limCond} we need the following lemma:
\begin{lemma}\label{lemma:fRegularity}
Let $w$ be an element of the Weyl group $W$.
\begin{enumerate}
\item[(a)] $\x^{\rho-w\rho}(1|w)(\x)$ is regular at the origin.
\item[(b)] $\x^{\rho-\sigma_iw\rho}\left(\smfrac{1}{x_i}|w\right)(\x)$ is
  regular at the origin for $i=1,2,\ldots, r$.
\end{enumerate}
\end{lemma}

\begin{proof}
The proof of the lemma is by induction on the length of $w$.  If
$w$ is the identity element, (a) and (b) above are trivial.
Suppose (a) and (b) are
true for $w_0\in W$ and that $i$ is such that $l(\sigma_iw_0)=l(w_0)+1$.
Then
\begin{eqnarray}\label{eqn:siw}
(1|\sigma_iw_0)(\x) &=& \left.\left(\frac{qx_i-1}{qx_i(1-x_i)}\right)
\right|w_0 \\
&=&\left.\left[g_1(x_i)+g_2(x_i)\smfrac{1}{x_i}\right]\right| w_0
\end{eqnarray}
where
$$g_1(x_i)=\frac{q-1}{q(1-x_i^2)} \quad \text{and}\quad
g_2(x_i)=\frac{qx_i^2-1}{q(1-x_i^2)} $$
are both even functions of $x_i$.  Therefore, by Lemma
\ref{lemma:Waction} (c),
$$(1|\sigma_iw_0)(\x) =
g_1(w_0\x)(1|w_0)(\x)+g_2(w_0\x)(\smfrac{1}{x_i}|w_0)(\x)
.$$
Since $g_i(w\x)$ is regular at the origin for $i=1,2$, to finish the
proof of (a) we must show that
$$\x^{\rho-\sigma_iw_0\rho}(1|w_0)(\x) \quad \text{and}\quad
\x^{\rho-\sigma_iw_0\rho}(\smfrac{1}{x_i}|w_0)(\x)$$
are both regular at the origin.  The second statement term is regular
by virtue of the inductive hypothesis.  As for the first, by induction
it suffices to show that
$\rho-\sigma_iw_0\rho\succeq\rho-w_0\rho$,
or equivalently, by Lemma \ref{lemma:Phi} (d), that
\begin{equation}\label{rootineq}
\sum_{\alpha\in \Phi(w_0^{-1})}\alpha \preceq
\sum_{\alpha\in \Phi(w_0^{-1})}\sigma_i(\alpha) +\alpha_i.
\end{equation}
In fact we claim that for each $\alpha\in \Phi(w_0^{-1})$ either
\begin{equation}\label{wialpha}
\alpha\preceq \sigma_i\alpha \mbox{\ or\ }  \sigma_i\alpha\in
\Phi(w_0^{-1}).
\end{equation}
Indeed, we know that
$\alpha-\sigma_i\alpha$ must be an integral multiple of $\alpha_i$, say
$\alpha-\sigma_i\alpha=n\alpha_i$.  If $n\leq 0$ then the first
alternative in \eqref{wialpha} holds.  If $n> 0, $ then
\begin{multline*}
w_0^{-1}\sigma_i(\alpha)=w_0^{-1}(\sigma_i\alpha-\alpha+\alpha)=\\
w_0^{-1}(\alpha)+w_0^{-1}(\sigma_i\alpha-\alpha)=w_0^{-1}\alpha-
nw_0^{-1}(\alpha_i)\in\Phi^-.
\end{multline*}
Now $\alpha$ is in $\Phi(w_0^{-1})$
and $\alpha_i$ is not.  Therefore $w_0^{-1}\sigma_i(\alpha)$ is in
$\Phi^-$ and $\sigma_i(\alpha)\in\Phi(w_0^{-1})$.
The proof of (b) is similar.
\end{proof}

\begin{proofof}{ of Theorem \ref{thm:invFunc}}
Let $f$ be defined as in \eqref{eqn:defoff}.  To complete the proof of
Theorem \ref{thm:invFunc}, we verify that $f$ satisfies the limiting
condition \eqref{limCond}.

Fix an index $i$ with neighbors $j_1, \ldots, j_k$.  Let $W_0$ be the
subgroup of $W$ generated by the $\sigma_j$ with $j\neq i$ and $j\neq
j_1,\ldots j_k$.  If we set $x_{j_1}=\cdots=x_{j_r}=0$ in
\begin{equation}\label{eqn:sum}
f(\x)=\Delta(\x)^{-1}\sum_{w\in W} j(w,\x) (1|w)(\x),
\end{equation}
then Lemmas \ref{lemma:Phi}(e), \ref{lemma:j}, and
\ref{lemma:fRegularity}(a)
imply that every summand in \eqref{eqn:sum} vanishes except for those
with $w$ in the group generated by $\sigma_i$ and $W_0$.
Since $\sigma_i$ is in the centralizer of $W_0$, \eqref{eqn:sum}
becomes
$$\Delta(\x)^{-1}\sum_{w\in W_0} [1+j(\sigma_i,\x) (1|\sigma_i)]|w.$$
The term in the brackets equals
$$1+\frac{qx_i^2-x_i}{1-x_i}=
\frac{1-qx_i^2}{1-x_i}$$
by \eqref{eqn:siw}.  Since each $w\in W_0$ is composed of reflections
$\sigma_j$ for $j$ neither neighboring nor equal to $i,$ the term
$\frac{1-qx_i^2}{1-x_i}$ can be pulled outside the summation, and leaves
behind a factor of $1/(1-x_i)$ after $1-qx_i^2$ cancels with the same
term in $\Delta$.  This completes the proof of Theorem
\ref{thm:invFunc}.
\end{proofof}

For use in the following sections, we establish some further properties
of the invariant function $f$.  Write $f(\x)=f(\x;q)$ as a power
series in the $x_i:$
\begin{equation}\label{eqn:psf}
f(\x;q)=\sum_{k_1,\ldots,k_r\geq 0} a(k_1,\ldots, k_r;q)x_1^{k_1}\cdots
x_r^{k_r} .
\end{equation}
We will often write $f(\x)$ or $a(k_1,\ldots,k_r)$ when the dependence
on $q$ is not relevant.  The main fact about the $q$-dependence
relevant for us is the following:

\begin{prp}\label{prp:polyBound}
For $\Phi$ fixed, there exists constants $C_1,C_2>0$
  such that $a(k_1, \dots,k_r;q)< C_1 q^{C_2|k|}$,
where $|k|:=k_1+\cdots+k_r$.
\end{prp}

\begin{proof}
From the definition of $f$, it is clear that its numerator is
polynomial in $q$ and that its denominator is a finite product of
terms of the form $(1-q^{l_0}x_1^{l_1}\cdots x_r^{l_r})$ for some
positive integers $l_i$.  
Expanding this out in a geometric series gives us the polynomial bound
in $q$.
\end{proof}

The reason for the introduction of the Weyl group action
\eqref{wiaction} and the relevance to $L$-functions
will be made more clear in the next section.  We conclude this section
by explaining a consequence of the $W$-invariance of
the function $f$.   Take the power series expansion of $f$ in $r-1$ of
the variables $x_i$.  Thus the coefficients of this expansion will be
functions of the one remaining variable, $x_{j_0}$, say.  The
invariance of $f$ under $\sigma_{j_0}$ will force these coefficients
to satisfy certain functional equations.  We make this explicit.

\begin{prp} \label{prp:feqnOfCPs}
Fix $q$ and an index $j_0$. Let
$$\hat k=(k_1,\ldots, k_{j_0-1}, k_{j_0+1}, \ldots,k_r)$$
be an $(r-1)$-tuple of nonnegative integers.  Define
$$T(x_{j_0};\hat k)=\sum_{k_{j_0}=0}^\infty a(k_1,\ldots, k_{j_0-1},k_{j_0},
k_{j_0+1},\ldots, k_r) x_{j_0}^{k_{j_0}}.$$
Let $n(\hat k)=\dsum{j\:adj(j,j_0)}{} k_j$.
\begin{enumerate}
\item[(a)]  If $n(\hat k)=2\gamma$ is even, then
$$(1-x)T(x_{j_0};\hat k)=(1-1/(qx_{j_0}))
(x_{j_0}\sqrt{q})^{2\gamma}T\left(\frac{1}{qx_{j_0}};\hat k\right).$$
\item[(b)]If $n(\hat k)=2\gamma+1$ is odd, then
$$T(x_{j_0};\hat k)=
(x_{j_0}\sqrt{q})^{2\gamma}T\left(\frac{1}{qx_{j_0}};\hat k\right).$$
\item[(c)]  Let $C_1,C_2$ be the constants of Proposition
  \ref{prp:polyBound}.    For $|x_{j_0}|<q^{-C_2}$, we have
$$|T(x_{j_0},\hat k)| < C_1 q^{-C_2|\hat k|},$$
where $|\hat k|=\dsum{j\neq j_0}{} k_j$.
\end{enumerate}
\end{prp}

\begin{proof}From \eqref{wiaction} and the invariance of $f$ under
  $\sigma_{j_0}$, we know that
  $(1-x_{j_0})f_{j_0}^+(\x)=(1-\smfrac{1}{qx_{j_0}}) )
f_{j_0}^+(\sigma_{j_0}\x)$.  Comparing the coefficients of
$$f_{j_0}^+(\x)=\sum_{\hat k\ :\ n(\hat k)\ \text{even}}T(x_{j_0};\hat k)
\prod_{j\neq j_0} x_j^{k_j}$$
and
$$ f_{j_0}^+(\sigma_{j_0}\x)=\sum_{\hat k\ :\ n(\hat k)\ \text{even}}
T\biggl(\frac{1}{qx_{j_0}};\hat k\biggr)
\biggl(\prod_{j\neq j_0} x_j^{k_j}\biggr)\biggl( \prod_{j:j,j_0 \adj}
(x_{j_0}\sqrt{q})^{k_j}\biggr)$$
yields (a).  The proof of (b) follows after a similar comparison of
$f_{j_0}^-(\x) $ and $f_{j_0}^-(\sigma_{j_0}\x)$.

\end{proof}

\begin{remark}\label{relation-to-wmd}
The rational functions of Theorem \ref{thm:invFunc} will be used to
define the $p$-parts of the multiple Dirichlet series of the following
section. (Here $p$ is a prime of norm $q.$)  An alternative description
of the $p$-parts of multiple Dirichlet series is given in the papers
\cite{wmd1, wmd2, wmd3}.  The first two of these papers deal with {\it
  stable} Weyl group multiple Dirichlet series constructed from
$n^{th}$ order characters and Gauss sums.  As noted in the
introduction, the series studied in this paper (the $n=2$ case) fall
outside the stable range provided $\Phi\neq A_2.$

To conclude this section, we describe the precise connection
between the $p$-part polynomial of \cite{wmd1, wmd2, wmd3} and
invariant rational function $f$ constructed above.  Our function
$f$ consists of both the $p$-part polynomial and the $p$-part of
the normalizing zeta factors of \cite{wmd1, wmd2, wmd3}.  In
Eq.~(30) of \cite{wmd2}, the normalizing zeta factor of the
quadratic multiple Dirichlet series associated to the root system
$\Phi$ of rank $r$ is defined to be
\begin{equation}\label{eqn:nzf}
\prod_{\alpha\in\Phi^+} \zeta(2\langle \alpha, {\bf
s}\rangle-d(\alpha)+1).
\end{equation}
Here, ${\bf s}$ is an $r$-tuple of complex numbers and
$$\langle\alpha,{\bf s}\rangle=\alpha_1 s_1+\cdots+\alpha_r s_r.$$
(Note: the formula of \cite{wmd2} is related to ours by the change
of variable $s_i\mapsto 2s_i-1/2.$) Thus, setting $x_i=q^{-s_i}, $
this product of zeta functions has $p$-part
\begin{equation}\label{eqn:denom}
D({\bf x})=\prod_{\alpha\in \Phi^+}
(1-q^{d(\alpha)-1}\x^{2\alpha}).
\end{equation}
Then $f(\x)D(\x)$ is a polynomial in the $x_i.$  After making
the change of variable $x_i\mapsto x_i\sqrt{q},$ this is the
$p$-part polynomial of \cite{wmd1, wmd2, wmd3}.

Let us compare our Examples \ref{ex:fa2} and \ref{ex:f} with the
formulas of \cite{wmd1, wmd2, wmd3}. We begin with the $A_2$
series, \eqref{eqn:a2series}.  Multiply $f_{A_2}(x,y)$ by
$(1-x^2)(1-y^2)(1-qx^2y^2).$ The result $N(x,y;A_2)$ is a sum of 6
terms which correspond to the 6 elements of the Weyl group $W.$
Make the change of variable $x\to x\sqrt{q}, y\to y\sqrt q$ in
$N(x,y;A_2)$ to get
$$ 1+\sqrt{q} x+\sqrt{q} y-q^{3/2}x^2y-q^{3/2}xy^2+q^2x^2y^2.$$
Then the coefficient of $x^{k_1}y^{k_2}$ is precisely the coefficient
$H(p^{k_1}, p^{k_2})$ given in (13) of \cite{wmd1}, after replacing
$g(1, p)$ by $\sqrt q$ and $g(p,p^2)$ by $-q.$  Thus, in this stable
example, our result is identical to the result of \cite{wmd1}.

Turning to Example \ref{ex:f},  multiplying  $f_{A_3}$ by
$$(1-x_1^2)(1-x_2^2)(1-x_3^2)(1-qx_1^2x_2^2)(1-qx_2^2x_3^2)
(1-q^2x_1^2x_2^2x_3^2)$$
yields a sum of 26 terms. After changes of variables as in the
paragraph above, 24 of these terms correspond to the 24 elements
of the Weyl group of $A_3$ under the association (6) of \cite
{wmd2}. However, (6) of \cite{wmd2} is intended to be applicable
only in the stable case; the missing 2 terms are a manifestation
of the instability of this example.

To investigate the connection between $f$ and the Weyl group, consider 
the rational function $f_{0} = f_{0}(\x ;q)$ from \eqref{defoff0}.
Expand $f_{0}$ as a power series in the variables $x_{i}$:
\begin{equation}\label{f0aspowerseries}
f_{0} (\x ;q) = \sum_{k_1,\ldots,k_r\geq 0} a(k_1,\ldots, k_r;q)x_1^{k_1}\cdots
x_r^{k_r}.
\end{equation}
It is not difficult to see that \eqref{f0aspowerseries} contains terms
in bijection with the Weyl group.  Indeed, consider the function
$f_{0} (\x ; 1)$ obtained by formally setting $q=1$ and applying the
definition \eqref{defoff0}.  If $q=1$, then the $W$-action
\eqref{wiaction} simplifies considerably, and one readily computes
\[
f_{0} (\x ;1) = \sum_{w\in W} (-1)^{l (w)+d (\rho -w\rho)}\x^{\rho -w\rho}.
\]
Since
$\rho$ lies in the interior of the Weyl chamber, it follows that the
monomials $\x^{\rho -w\rho}$ are all distinct.  This proves that
\eqref{f0aspowerseries} contains terms in bijection with $W$.  This
also shows that, as functions in $q$, the coefficents of the unstable
terms vanish when $q=1$.

Finally, we note that Brubaker, Bump, Friedberg and Hoffstein
\cite{wmd3} have given a conjecture for the $p$-parts, applicable for
all $n$ when $\Phi=A_r.$ In this conjecture the terms of the numerator
are parametrized not by Weyl group elements, but rather by
Gelfand--Tsetlin patterns of rank $r$ with top row $(r, r-1, \ldots,
2,1)$.  These terms include terms parametrized by the Weyl group; as
monomials they coincide with the $\x^{\rho -w\rho}$ from above.
Moreover, the additional unstable terms in their conjectural $p$-parts
satisfy a remarkable geometric property.  Let $P$ be the convex
polytope obtained by taking the convex hull of the points $\rho
-w\rho, w\in W$ in the vector space $\Lambda_{\Phi}\otimes \R$.  Then
the unstable terms are supported on monomials $\x ^\alpha $ with
$\alpha \in \Lambda_{\Phi}$ and lying in $P$.

The authors of \cite{wmd3} provide much convincing evidence for their
conjecture, including verification that for $n=2,$ the conjecture
agrees with our results for $A_r, r\leq 5.$ Unfortunately, our methods
do not readily provide a means to attack their conjecture as it
appears difficult to extract the coefficients of the numerator of the
rational function of Theorem \ref{thm:invFunc} from the definition
\eqref{eqn:defoff}, and because lots of cancellation occurs during the
averaging process.  The connection between our construction and that
of \cite{wmd3} is currently under investigation by the authors in
joint work with Bump and Friedberg.
\end{remark}

\section{Definition of the quadratic Weyl group multiple Dirichlet
series}
\label{sec:mds}

We continue to let $\Phi$ denote an irreducible simply-laced root
system of rank $r$.   We recall our convention on the ordering of the
indices: for each $j$, the $i$ which are adjacent to $j$ are
either all less than $j$ or all greater than $j$.

Let
$$\Psi=(\psi_1,\psi_2,\ldots,\psi_r)$$
be  a collection of $r$ id\`ele class characters unramified outside of $S$.
Given a collection ${\bf I}=(I_1,\ldots,I_r)$ of ideals in $\I(S)$ we denote
by $\Psi({\bf I})$ the product
$$\prod_{i}\psi_i(C_i).$$
and by $H({\bf I})$ the coefficient $H(I_1,I_2,\ldots,I_r)$
defined below.

\begin{defin} \label{def:H} The coefficient $H(I_1,I_2\ldots,I_r)$ is
  defined by the following two conditions:
\begin{enumerate}
\item Suppose $\mathbf{I} = (P^{k_{1}},\dotsc , P^{k_{r}})$, where $P$
is a fixed prime ideal of norm $q$.  Then
\begin{equation*}H(P^{k_1},\ldots,P^{k_r})=a(k_1,\ldots, k_r;q).
\end{equation*}
\item  Given ideals $I_j, I_j^\prime\in \I(S)$ with $(I_1I_2\cdots
  I_r, I_1^\prime I_2^\prime\cdots I_r^\prime)=1$ we have
\begin{equation*}
\frac{H(I_1I_1^\prime, \ldots, I_rI_r^\prime )}
{H(I_1, \ldots, I_r)H(I_1^\prime, \ldots, I_r^\prime )}=
\prod_{\substack{i,j \adj\\i<j}} \prs{I_i}{I_j^\prime}
\prs{I_i^\prime} {I_j}
\end{equation*}
\end{enumerate}
\end{defin}

Note that the second condition and Proposition \ref{prp:polyBound}
imply the bound
\begin{equation}\label{eqn:HBound}
|H(I_1,\ldots, I_r)|<\!\!< |I_1\cdots I_r|^C
\end{equation}
for some constant $C$.
If the ideals $I_1,\ldots, I_r$ is are pairwise
relatively prime, then $H({\bf I})$ has an especially simple form:

\begin{lemma}
If the ideals $I_1,\ldots , I_r\in \I(s)$ are pairwise relatively prime,
then
$$H(I_1, \ldots, I_r)=\prod_{\substack{i,j \adj\\i<j}} \prs{I_i}{I_j}.$$
\end{lemma}

\begin{proof}
We have
\begin{equation*}
H(I_1,I_2,\ldots, I_r) = H( I_1,1,\ldots,1)
H( 1,I_2,\ldots, I_r)
 \prod_{\substack{i,1 \adj}} \prs{C_1}{C_i}.
\end{equation*}
Now use the fact that $H(C,1\ldots,1)=1$ and induct.
\end{proof}

We may finally define the family of multiple Dirichlet series that is the
main subject of this paper.  For an $r$-tuple $\s=(s_1,\ldots,s_r)$ of
complex numbers, define
\begin{equation}\label{def:ZS}
Z_S(\s,\Psi)=\sum_{\II=(I_1,\ldots, I_r)\in \mathcal I (S)^r}
\frac{\Psi(\II)H(\II)}
{\prod_j |I_j|^{s_j}}
\end{equation}
By the \eqref{eqn:HBound} we see that the sum defining $Z_S(\s,\Psi)$
will converge absolutely for $\re(s_j)$ sufficiently large, $1\leq
j\leq r$.

We will find it convenient to extend this definition to allow linear
combinations of id\`ele class characters in place of $\Psi$.  If
$$\Xi=\sum b_\Psi \Psi$$
for some collection of complex numbers $b_\Psi$, we define
$$Z_S(\s, \Xi)=\sum b_\Psi Z_S(\s, \Psi).$$
In the particular applications we have in mind, the $r$-tuple $\Xi$
will consist of combinations of id\`ele class characters and
characteristic functions $\delta_E$ for classes $E$  in $R_C$.

\begin{remark}
The coefficient function $H$ is similar to but slightly different from
the function of the same name in \cite{wmd1,wmd2,wmd3}.  To compare
the two, denote the function in \cite{wmd1,wmd2,wmd3} by $H_{\WMD}$.  As
explained in Remark \ref{relation-to-wmd}, the coefficient generating
function $f(\x;q)$ contains both the $P$-part polynomial of
\cite{wmd1, wmd2, wmd3} and the normalizing zeta factor
\eqref{eqn:nzf}.  Therefore, we expect the equality
\begin{multline}\label{eq:conjecturedequality}
 D(\x)\sum_{k_1,\ldots,k_r} H(P^{k_1},\ldots,P^{k_r}) x_1^{k_1}\cdots
x_r^{k_r}\\ = \sum_{k_1,\ldots,k_r} H_{\WMD}(P^{k_1},\ldots,P^{k_r})
y_1^{k_1}\cdots y_r^{k_r}
\end{multline} 
where $x_i\sqrt{q}=y_i$ and $D$ is the denominator given in
\eqref{eqn:denom}.  The coefficients on the right hand side are to be
understood to mean those defined in \cite{wmd1} when $\Phi=A_2$ and to
mean those conjectured in \cite{wmd3} when $\Phi=A_r$ for $r\geq 3,
n=2$.  As mentioned in Remark \ref{relation-to-wmd}, we have checked
equality of \eqref{eq:conjecturedequality} for $r\leq 5$, $n=2$.
\end{remark}

\section{Functional equations and analytic continuation}
\label{sec:fes}

In this section we show that the family of multiple Dirichlet series
$Z_S(\s,\Psi)$ as $\Psi$ ranges over $r$-tuples of quadratic id\`ele class
characters unramified outside of $S$ satisfies a group of functional
equations isomorphic to $W$, the Weyl group of the root system
$\Phi$.  Summing over the $j_0^{th}$ index in the series \eqref{def:ZS}
defining $Z_S(\s, \Psi)$ will produce an $L$-function having a
functional equation as $s_{j_0}\mapsto 1-s_{j_0}$. This functional equation
will induce a functional equation in the multiple Dirichlet
series relating the values at $\s=(s_1,\ldots, s_r)$ to the values at
$\sigma_{j_0}\s=(s_1^\prime,\ldots,s_r^\prime)$, where
\begin{equation}\label{eqn:actionOnS}
s_j^\prime=\left\{\begin{array}{ll}
s_j+s_{j_0}-1/2 & \mbox{\ if $j$ and $j_0$ are adjacent,}\\
1-s_{j_0} & \mbox{\ if $j=j_0$, and}\\
s_j & \mbox{\ otherwise.}
\end{array}\right.
\end{equation}
 These functional equations are involutions generating the
group of functional equations of $Z_S(\s, \Psi)$.  Note that if we set
$x_j=q^{-s_j}$, then this action corresponds to the action
\eqref{eqn:wiaction1}  of $W$ on  $\x=(x_1,\ldots, x_r)$
by the variable change $x_j=q^{-s_j}$.

We now exhibit the functional equations in detail.  Fix an index
$j_0$.   Then summing \eqref{def:ZS} over this index first produces
\begin{equation}\label{eqn:sumOverOneIndex}
\sum_{\substack{j=1,\ldots, r\\j\neq j_0}}
\sum_{I_j\in \mathcal I (S)}
\frac{\prod_{j\neq j_0}\psi_j(I_j)}{ \prod_{j\neq j_0}
|I_j|^{s_j}}
\cdot \sum_{I_{j_0}\in \I(S)} \frac{H(I_1,\ldots,I_{j_0-1},I_{j_0},
  I_{j_0+1},\ldots, I_r)}{|I_{j_0}|^{s_{j_0}}}\psi_{j_0}(I_{j_0}).
\end{equation}
Our goal is to express the innermost sum as the product of
a partial $L$-series with a Dirichlet polynomial, and to exhibit the
precise functional equation that it satisfies.

Let $N=\prod_{j\neq j_0} I_j$ and let $M=\prod_{j:j, j_0 \adj} I_j$.
We will assume that $j_0>j$ for all
indices $j$ adjacent to $j_0$. Setting
$\psi=\psi_{j_0}$ and $s=s_{j_0}$,
we begin by removing the ideals relatively prime to
$N$ from the inner sum above:
\begin{eqnarray*}
\lefteqn{\sum_{I_{j_0}\in\I(S)} \frac{H(I_1,\ldots,I_{j_0},
  \ldots, I_r)}{|I_{j_0}|^s}\psi(I_{j_0})} \\
&=&\sum_{I|N^\infty}\sum_{(J,N)=1 }\frac{H(I_1,\ldots,IJ,\ldots, I_r)}
{|IJ|^s}\psi(IJ)\\
&=&\sum_{I|N^\infty}\frac{H(I_1,\ldots,I,\ldots, I_r)}{|I|^s}\psi(I)
\left[
\sum_{(J,N)=1 }\frac{\psi(J)}{|J|^s}\prod_{j:j,j_0 \adj} \prs{I_j}{J}
\right]\\
&=& L_{S_N}(s,\psi\chi_M)
\sum_{I|N^\infty}\frac{H(I_1,\ldots,I,\ldots,I_r)}{|I|^s}
\psi(I),
\end{eqnarray*}
where $S_N$ is the set of places in $S$ together with the places
dividing $N$.  (If we had chosen $j_0$ such that all $j$ adjacent to
$j_0$ had been greater than $j_0$ the only difference would be that
the partial $L$-function in front would instead be associated to the
character $\psi\psi_M\chi_M$ where $\psi_M$ is the (unramified outside
$S$) id\`ele class character given by $J\mapsto  \alpha(I,J)$ which
depends only on the class of $M$ in $R_C$.)
The sum over $I|N^\infty$ decomposes as a product over the primes
dividing $N$.  Let $P$ be a prime divisor of $N$. Let $\beta_j$ be
the order of $P$ in $I_j$ and let $I_j^{(P)}$ denote the part of $I_j$
relatively prime to $P$.  Thus $I_j=I_j^{(P)}P^{\beta_j}$.  Then
\begin{eqnarray}\label{eqn:Hsum1}
\lefteqn{\sum_{I|N^\infty}\frac{H(I_1,\ldots,I,\ldots,I_r)}{|I|^s}
\psi(I)  }\\
&=& \sum_{\substack{I|N^\infty\\(I,P)=1}}\sum_{k=0}^\infty
\frac
{H(I^{(P)}_1P^{\beta_1}, \ldots,IP^k, \ldots,I^{(P)}_rP^{\beta_r})}
{|I|^s|P|^{ks}}\psi(IP^k).\nn
\end{eqnarray}
Using the twisted multiplicativity,  the term
$H(I^{(P)}_1P^{\beta_1}, \ldots,IP^k, \ldots,I^{(P)}_rP^{\beta_r})$
in the numerator  can be pulled apart to yield
\begin{multline}
H(I^{(P)}_1, \ldots,I, \ldots,I^{(P)}_r)
H(P^{\beta_1}, \ldots,P^k, \ldots,P^{\beta_r})\\
\times \left[\prod_{\substack{i<j, \adj \\i,j\neq j_0}}
\prs{I_i^{(P)}}{P^{\beta_j}}\prs{P^{\beta_i}}{I_j^{(P)}}\right]
\left[\prod_{j:j,j_0 \adj} \prs{I_j^{(P)}}{P^k}\prs{P^{\beta_j}}{I}
\right].
\end{multline}
The first bracketed product of characters is a constant which can be
pulled outside the summation and will be ignored.  Summing over $k$ we
get
$$\sum_{k=0}^\infty \frac{H(P^{\beta_1}, \ldots,P^k,
  \ldots,P^{\beta_r})} {|P|^{ks}}\psi(P^k)\prod_{j:j,j_0 \adj}
  \prs{I_j^{(P)}}{P^k}
.$$
Thus, up to a constant of absolute value 1, \eqref{eqn:Hsum1} is
\begin{equation}\label{eqn:Hproduct}
\prod_{\substack{P|N\\P^{\beta_j}||I_j}}
\sum_{k=0}^\infty \frac{H(P^{\beta_1}, \ldots,P^k,
  \ldots,P^{\beta_r})} {|P|^{ks}}\psi(P^k)\prs{M^{(P)}}{P^k}.
\end{equation}

Write $M=M_0M_1^2M_2^2$ with $M_0$ squarefree and $(M_0M_1, M_2)=1$.
Therefore, $M_1$ consists of primes which divide $M$ to odd power and
$M_2$ of primes dividing $M$ to even order.  In further evaluating the
product \eqref{eqn:Hproduct}, we distinguish three cases:  $P$
relatively prime to $M$, $P$ divides $M$ to odd order, and $P$ divides
$M$ to even order.

{\em Case 1:  $P$ relatively prime to $M$.}  This means that for all
the neighbors $j$ of $j_0, P^{\beta_j}=1$.  By the limiting condition
of Theorem \ref{thm:invFunc} we conclude that
$$\sum_{k=0}^\infty \frac{H(P^{\beta_1}, \ldots,P^k,
  \ldots,P^{\beta_r})} {|P|^{ks}}\psi(P^k)\chi_M(P^k)$$
is a constant (independent of $s$) multiple of
$$(1-\psi(P)\chi_M(P)|P|^{-s})^{-1},$$
the $P$-part of the $L$-function $L(s, \psi\chi_M)$. The constant is
given by the $k=0$ term:
\begin{equation}\label{eqn:caseZeroConstant}
H(P^{\beta_1}, \ldots,P^k,
\ldots,P^{\beta_r})<<|P|^{C(\beta_1+\cdots+\beta_r)}.
\end{equation}

{\em Case 2: $P$ divides $M$ to odd order.}   Let the order of $P$ in
$M$ be $2\gamma+1$.  Let $\epsilon=\pm 1$ be $\psi(P)\prs{M^{(P)}}{P}$.
Thus the $P$-part of \eqref{eqn:Hproduct} is
$$H_P(s):=\sum_{k=0}^\infty \frac{H(P^{\beta_1}, \ldots,P^k,
  \ldots,P^{\beta_r})} {|P|^{ks}} \epsilon^k.$$
By virtue of the functional equation satisfied by $f_i^-$
(Proposition  \ref{prp:feqnOfCPs}), $H_P(s)$ satisfies
$$H_P(s)=|P|^{\gamma(1-2s)}H_P(1-s).$$  Taking the product over all
  $P$ dividing $M$ to odd order, we have
\begin{equation}\label{eqn:feHpOdd}
\prod_{\ord_P(M) odd} H_P(s) = |M_1|^{1-2s} \prod_{\ord_P(M) odd}
H_P(1-s).
\end{equation}

{\em Case 3: $P$ divides $M$ to even order.}
Let the order of $P$ in $M$ be $2\gamma$.  In this case,
$\chi_{M^{(P)}}=\chi_M$ since $\chi_M$ depends only on the squarefree
  part of $M$.  The  $P$-part of \eqref{eqn:Hproduct} is
$$H_P(s):=\sum_{k=0}^\infty \frac{H(P^{\beta_1}, \ldots,P^k,
  \ldots,P^{\beta_r})} {|P|^{ks}} \psi\chi_M(P^k).$$
Again by Proposition \ref{prp:feqnOfCPs}, this can be written as
$$(1-\psi\chi_M(P)|P|^{-s})^{-1} H_P(s),$$
where $H_P(S)$ satisfies
$$H_P(s)=|P|^{2\gamma(1/2-s)}H_P(1-s).$$
Taking the product over all $P$ dividing $M$ to even order, we have
\begin{equation}\label{eqn:feHpEven}
\prod_{\ord_P(M) even} H_P(s) = |M_2|^{1-2s} \prod_{\ord_P(M) even}
H_P(1-s).
\end{equation}

Putting together the 3 cases above, we get an expression for
\eqref{eqn:sumOverOneIndex} in terms of an $L$-function.

\begin{prp}\label{prp:LFunction}Fix ideals $I_j\in \I(S)$ for $j\neq
  j_0$.
Then
$$\sum_{I_{j_0}\in\I(S)} \frac{H(I_1,\ldots,I_{j_0},
  \ldots, I_r)}{|I_{j_0}|^s}\psi(I_{j_0})=
L^S(s,\psi\chi_M)Q(s)$$
where $Q(s)$ is a finite Euler product depending on the ideals $I_1,
\dots I_{j_0-1},I_{j_0+1},\ldots , I_r$ and the character $\psi$ which
  satisfies
\begin{equation}\label{eqn:feQ}
Q(s)=|M_1M_2|^{1-2s}Q(1-s).
\end{equation}
For $s>C_2$, there exists $C_3$ such that
$$|Q(s)| < |N|^{C_3}.$$
Here, $C_2$ is the constant from Proposition
\ref{prp:feqnOfCPs}.
\end{prp}

\begin{proof}
The only unproven part of the Proposition is the claim about the size
of $Q(s)$.  It follows from  Proposition \ref{prp:feqnOfCPs} that
 $$|Q(s)| < C_1^{\omega(N)}|N|^{C_2}$$
where $\omega(N)$ is the number of prime divisors of
$N$.  Hence we may take $C_3=C_2+\log_2 C_1$.
\end{proof}

Let
$$\hat L^S(s,\psi\chi_M):=\sum_{I_{j_0}\in\I(S)} \frac{H(I_1,\ldots,I_{j_0},
  \ldots, I_r)}{|I_{j_0}|^s}\psi(I_{j_0}).$$
Note that \eqref{eqn:feQ} forces $Q(s)$ to be a Dirichlet polynomial.
Therefore, $Q(s)$ is an entire function of $s$.  This implies that
$\hat L^S(s,\psi\chi_M)$,  has an analytic continuation to $s\in
\C$, with at most a simple pole at $s=1$.  This simple pole will exist
if and only if $\psi\chi_M$ is the trivial character.  Moreover, $\hat
  L^S(s,\psi\chi_M)$  will satisfy a functional equation as $s\mapsto
  1-s$.

\begin{prp}\label{prp:feLHat}There is a factor $A(s,\psi,E)$ depending only
  on $\psi$ and the class $E\in\E$ of $M$ such that
$$L_S(s,\psi\chi_M)\hat L^S(s,\psi\chi_M)=A(s,\psi,E) M^{1/2-s}
  L_S(1-s,\psi\chi_M)\hat L^S(1-s,\psi\chi_M). $$
In fact, $L_S(s,\psi\chi_M)$ also depends only on $\psi$ and the class
  of $M$ in $\E$.   The function $A(s,\psi,E)$ is of the form
  $A_0^{1/2-s}$ where $A_0=\frac{|\f_M|}{|M_0\f_E|}$.
\end{prp}

This is immediate from \eqref{eqn:feQ}, the functional equation
\eqref{eq:LFE} of the  $L$-function $L(s,\psi\chi_M)$ and the
description of the epsilon factor in Proposition
\ref{prop:epsilonFE}.  We emphasize that the fact that
$L_v(s,\psi\chi_M)$ depends only on $\psi$ and the class of $M$ in
$\E$ is true for both the archimedean and nonarchimedean places
$v\in S$---see the remark after Propostion \ref{prop:epsilonFE}.

In the usual way, we can use the functional equation of the preceding
proposition to obtain a convexity estimate for $\hat L^S:$
$$\hat L^S(s,\psi\chi_M)<< |N|^{C_3} |M|^{2C_2+1}$$
for $\re(s)>-C_2$, with the implicit constant depending on the set $S$
and $\im(s)$.  This final estimate allows us to analytically continue
$Z_S(\s,\Psi)$ slightly beyond the initial domain of absolute
convergence.

\begin{prp}\label{prp:ac}
 For each $j_0$, the multiple Dirichlet series $Z_S(\s,\Psi)$ has an
  analytic continuation to the the domain
$$\Omega=
\{(s_1,\ldots,s_r)\in \C^r: \re(s_{j_0})>-C_2, \re(s_j)>C_3+2C_2+2,
  \mbox{\ for \ } j\neq j_0\}.$$
\end{prp}

The actual constants $C_1, C_2, C_3$ are unimportant.  The point is
that the base of the  tube domain described in the previous proposition
is the
complement of a compact subset of the base of the orthant
$$X=\{\re(s_j)>0\mbox{\ for\ } j=1,2,\ldots, r\}.$$

Let $E$ be a class in $\E$. Let $\delta_{j_0,E}$ be the function on
$\I(S)^r$ defined by
$$\delta_{j_0,E}(I_1,\ldots, I_r)=
\left\{\begin{array}{ll}
1 & \mbox{if\ }\dprod{j:j,j_0 \adj}{}I_j\sim E \\
0 & \mbox{otherwise.}
\end{array}\right.
$$
Write
$$Z_S(\s,\Psi)=\sum_{E\in \E} Z_S(s,\Psi\delta_{j_0, E}).$$
Then
$$ Z_S(\s,\Psi\delta_{j_0, E})=\sum_{\substack{j=1,\ldots, r\\j\neq j_0}}
\sum_{I_j\in \mathcal I (S)}
\frac{\prod_{j\neq j_0}\psi_j(I_j)}{ \prod_{j\neq j_0}|I_j|^{s_j}}
\hat L^S(s,\psi\chi_M).$$
Write $B(s,\psi,\E)=L_S(s,\psi,\chi_M)$ and multiple $
Z_S(\s,\Psi\delta_{j_0, E})$ by this factor.  Then, using Proposition
\ref{prp:feLHat},  we have the
functional equation
$$ B(s,\psi,E)Z_S(\s,\Psi\delta_{j_0,
  E})=A(s,\psi,E)B(1-s,\psi,E)Z_S(\sigma_{j_0}\s,\Psi\delta_{j_0, E}).$$

Recall that the action of $W$ on $\s$ was given in \eqref{eqn:actionOnS}.
Summing over $E$ we get the functional equation for $Z_S(\s,\Psi)$.

\begin{theorem}\label{thm:feOfZ}  For each $j_0=1,2,\ldots, r$,
$$Z_S(\s,\Psi)=\sum_{E\in\E}A(s,\psi,E)\frac{B(1-s,\psi,E)}{B(s,\psi,E)}
Z_S(\sigma_{j_0}\s,\Psi\delta_{j_0, E}).$$
\end{theorem}

Let $\vec Z_S(\s)$ be the vector consisting of the $Z_S(\s,\Psi)$ as
$\Psi$ ranges over  $r$-tuples of quadratic id\`ele class characters
unramified outside of $S$.   Writing an arbitrary element $w\in W$ in
terms of the simple reflections, we may express Theorem
\ref{thm:feOfZ} as
\begin{equation}\label{eqn:feOfZ}
\vec Z_S(\s)=\Phi(\s;w) \vec Z_S(w\s)
\end{equation}
for some matrix $\Phi(\s;w)$.

\begin{theorem}\label{thm:main}
The function  $Z_S(\s,\Psi)$ has an analytic continuation to $\C^r$.
The collection of these functions as $\Psi$ ranges over
$r$-tuples of quadratic id\`ele class characters unramified outside of
$S$ satisfies a group of functional equations isomorphic to $W$.
This action of $W$ is given by Theorem \ref{thm:feOfZ} and
\eqref{eqn:feOfZ}. Finally, $Z_S(\s,\Psi)$ is analytic outside the
hyperplanes  $(w\s)_j=1$, for $w\in W, 1\leq j\leq r$. Here  $(w\s)_j$
denotes the $j^{th}$ component of $w\s$.
\end{theorem}

\begin{proof}
The argument is identical to that given in the proof of Theorem 5.9 of
\cite{wmd2}, and we do not repeat the details here.  However, for the
convenience of the reader, we give a sketch.  Using the functional
equations \eqref{eqn:feOfZ}, we may extend the domain of analyticity
of $Z_S(\s,\Psi)$ to translates of $\Omega$ by the group $W$.  The
union of the translates forms a tube domain in $\C^r$ whose base is
the complement of a compact subset of $\R^r$.  We may then apply
Bochner's theorem \cite{bo} to extend $Z_S(\s,\Psi)$ to an analytic
function on  all of $\C_r$.
\end{proof}

\end{document}